\documentclass[reqno, 10pt, centertags,draft]{amsart}
\usepackage{amsmath,amsthm,amscd,amssymb,latexsym,upref}
\usepackage[mathscr]{eucal}

\makeatletter
\def\theequation{\@arabic\c@equation}

\newcommand{\bbN}{{\mathbb{N}}}
\newcommand{\bbR}{{\mathbb{R}}}

\newcommand{\bbC}{{\mathbb{C}}}

\newcommand{\R}{\mathbb{R}}

\newcommand{\cB}{{\mathcal B}}

\newcommand{\cH}{{\mathcal H}}

\newcommand{\lb}{\label}
\newcommand{\f}{\frac}

\newcommand{\ol}{\overline}

\newcommand{\loc}{\text{\rm{loc}}}

\newcommand{\ran}{\text{\rm{ran}}}

\newcommand{\dom}{\text{\rm{dom}}}

\newcommand{\slimes}{\text{\rm{l.i.m.}}}

\newcommand{\bi}{\bibitem}

\renewcommand{\Im}{\text{\rm Im}}

\newtheorem{theorem}{Theorem}

\theoremstyle{definition}

\begin{document}

\title[Hill Operators and Spectral Operators of Scalar Type]{When Is A 
Non-Self-Adjoint Hill Operator \\ A
Spectral Operator Of Scalar Type?}
\author[F.\ Gesztesy, and V.\ Tkachenko]{Fritz
Gesztesy and Vadim Tkachenko}
\address{Department of Mathematics,
University of Missouri, Columbia, MO 65211, USA}
\email{fritz@math.missouri.edu}
\urladdr{http://www.math.missouri.edu/personnel/faculty/gesztesyf.html}
\address{Department of Mathematics,
Ben Gurion University of the Negev, Beer--Sheva 84105, Israel}
\email{tkachenk@cs.bgu.ac.il}
\thanks{Based upon work supported by the US National Science
Foundation under Grant No.\ DMS-0405526 and the Research Council and
the Office of Research  of the University of Missouri--Columbia.}
\thanks{Scientific Field: Ordinary Differential Equations.}

\begin{abstract}

We derive necessary and sufficient conditions for a one-dimensional periodic
Schr\"odinger (i.e., Hill) operator $H=-d^2/dx^2+V$ in $L^2(\bbR)$ to be a
spectral operator of scalar type. The conditions demonstrate the remarkable
fact that the property of a Hill operator being a spectral operator is
independent of smoothness (or even analyticity) properties of the potential
$V$. 
\end{abstract}

\maketitle

{\bf R\'ESUM\'E. Quand un op\'erateur de Hille non-autoadjoint est-il un
operateur sp\'ectral de type scalaire?} Nous derivon des conditions
n\'ecessaires et suffisantes pour qur  l'op\'e\-ra\-te\-ur de
Schr\"odinger (i.e., l'op\'erateur de Hill) 
$H=-d^2/dx^2+V$ 
dans $L^2(\R)$ soit un op\'erateur sp\'ectral de type scalaire.
Les conditions  montrent que cette propri\'et\'es ne d\'epend pas des
propri\'et\'es  diff\'erentieles (ou analytiques) du potentiel $V$.

\bigskip


We consider (maximally defined) Hill operators
$H=- \frac{d^2}{dx^2} + V$, $x\in\mathbb{R}$, 
in $L^2(\bbR)$ with complex-valued $\pi$-periodic potentials $V$ such
that $V\in L^2([0,\pi])$.

Then two fundamental problems, central to the spectral theory of such
operators, are the following:

$\mathbf{(I)}$ What is the spectrum, $\sigma(H)$, of $H$, in particular,
what is the nature and geometry of $\sigma(H)$\,?

$\mathbf{(II)}$ What kind of spectral decomposition (or resolution of the
identity), if any, is generated by $H$ in the space $L^2(\mathbb{R})$\,?

Both problems are completely solved for the case of real-valued potentials,
that is, for self-adjoint Hill operators. As to the first problem, the
Floquet--Bloch theory states that for every such operator there exists
a sequence of real numbers
$$
\lambda_0^+< \lambda_1^-\leq\lambda_1^+< \cdots <\lambda_k^-
\leq\lambda_k^+< \cdots
$$
such that the spectrum of $H$ is purely absolutely continuous and has the
form
\begin{equation}\label{2}
\sigma(H)=\bigcup\limits_{k=0}^\infty \; [\lambda_k^+,\lambda_{k+1}^-].
\end{equation}

A solution of the second problem for self-adjoint Hill operators
was provided by Titchmarsh \cite{Ti50} around 1950. He  obtained an
explicit formula for the resolution of the identity generated by $H$ in
terms of the fundamental matrix
$$
U(z,x)=
\begin{pmatrix}
\theta(z,x)   & \phi(z,x)\\
\theta'(z,x)   & \phi'(z,x)
\end{pmatrix},
\quad U(z,0)=I_2, \; z\in\bbC,
$$
of solutions of the differential equation
\begin{equation}
-y''(z,x) +V(x) y(z,x)=zy(z,x),  \label{20}
\end{equation}
where prime $\prime$ denotes the derivative with respect to $x\in\bbR$.

According to Titchmarsh \cite{Ti50}, for every element $f\in
L^2(\mathbb{R})$, the orthogonal $L^2$-representation
\begin{equation}\label{4}
f(x)=\sum\limits_{k=0}^\infty \f{(-1)^k}{2\pi}\int_{
\lambda_k^+}^{\lambda_{k+1}^-} \frac{d\lambda}
{\sqrt{1-\Delta_+(\lambda)^2}} \, \Phi(\lambda,x;f)
\end{equation}
holds with 
\begin{align}
\begin{split}
\Phi(\lambda,x;f)&=\phi(\lambda,\pi)\theta(\lambda,x)F_1(\lambda)
-\theta'(\lambda,\pi)\phi(\lambda,x)F_2(\lambda)  \lb{1.9} \\
& \quad -\Delta_-(\lambda)\theta(\lambda,x)F_2(\lambda)
-\Delta_-(\lambda)\phi(\lambda,x)F_1(\lambda), \quad \lambda\in\sigma(H),
\end{split}
\end{align}
where $\Delta_\pm(z)=[\theta(z,\pi)\pm\phi'(z,\pi)]/2,\; z\in\bbC$, 
$\sqrt{1-\Delta_+(\lambda)^2}\geq 0$, and $F_1(\lambda)=\int_\mathbb{R}
dy\,f(y) \theta(\lambda,y)$, $F_2(\lambda)=\int_\mathbb{R} dy\,f(y)
\phi(\lambda,y)$, for $\lambda\in\sigma(H)$.
As a corollary, for every closed Borel set
$\sigma\subset\sigma(H)$  the operator
$P(\sigma)$ defined by
\begin{equation}\label{5}
(P(\sigma)f)(x)=\frac{1}{2\pi}\int_{\sigma} \f{d\lambda}
{\sqrt{1-\Delta_+(\lambda)^2}} \,\Phi(\lambda,x;f)
\end{equation}
is an orthogonal projection of $L^2(\mathbb{R})$ onto the closed subspace
$\ran(P(\sigma))$ invariant with respect to $H$, such that the spectrum of
its restriction to $\ran(P(\sigma))$ is contained in $\ol\sigma$.

While the case of self-adjoint Hill operators is under
complete control, the case of non-self-adjoint Hill operators
offers a remarkable complexity. The first result in this context 
obtained by Serov \cite{Se60} reads:

\smallskip
\noindent
{\em The spectrum of a Hill operator with a complex-valued
potential $q\in L^2([0,\pi])$ 
coincides with the set }
$\sigma(H)=\{\lambda\in\mathbb{C}\,|\, \Delta_+(\lambda)\in[-1,1]\}$.
%
According to this result 
the spectrum $\sigma(H)$ of $H$ is
formed by a system of analytic arcs in the complex plane that may
intersect in inner points as shown in \cite{PT91a} (see also \cite{GW95}).
For necessary and sufficient conditions on a set $\Sigma\subset\bbC$ to be
the spectrum of some periodic Schr\"odinger operator $H$ with periodic
potential $q\in L^2_{\loc}(\bbR)$ in terms of a certain class of
Riemann surfaces, we refer to \cite{Tk96}.

While the spectrum of non-self-adjoint Hill operators has been understood
since 1960, much less was known about the spectral decompositions
generated by such non-self-adjoint operators. At first this appears to be
unusual since all ingredients of the Titchmarsh formula \eqref{4} are
also present in the non-self-adjoint situation and hence it would seem
natural to use them for a corresponding spectral analysis of the
non-self-adjoint case. Such attempts, however, immediately meet essential
obstacles. 

According to general principles of spectral analysis, to adjust the formula
\eqref{4} to the non-self-adjoint case one replaces the intervals 
$[\lambda_k^+,\lambda_{k+1}^-]$ by
the spectral arcs comprising $\sigma(H)$ and instantly meets the first
difficulty at the points of intersection of such arcs. 
If $\lambda_0$ is such a point  and $^\bullet$ denotes the derivative with
respect to $z$, then a change of variables from $\lambda\in[\lambda_k^+,
\lambda_{k+1}^-]$ (now denoting a spectral arc) to $t\in [0,2\pi]$, using
$\Delta_+(\lambda)=\cos (t)$,  leads to integrals of the type $\int_{T}{dt}
(\Delta_+^\bullet(\lambda(t)))^{-1} \, \Phi(\lambda(t),x;f)$, 
where $T\subseteq[0,2\pi]$ is an interval containing $t_0$.
Convergence of these integrals for
arbitrary $f\in L^2([0,\pi])$ and their operator properties depend on
the behavior of the meromorphic functions
\begin{equation}\label{111}
\frac{\phi(z,\pi)}{\Delta_+^{\bullet}(z)},\quad
\frac{\theta'(z,\pi)}{\Delta_+^{\bullet}(z)},\quad 
\frac{\Delta_-(z,\pi)}{\Delta_+^{\bullet}(z)}
\end{equation}
participating in \eqref{1.9}.

Clearly, Hill operators $H$ commute with the
operator of translation
by $\pi$. McGarvey \cite{Mc62} initiated the study of general
operators of such a type within the framework of the theory of spectral
operators in the sense of  Dunford \cite{DS88}. A
closed operator $T$ with domain $\dom(T)\subseteq\mathcal{H}$ is called a
{\em spectral operator} if there exists a countably additive
projection-valued  measure $E_T(\cdot)$ defined on the Borel subsets
$\cB$ of $\mathbb{C}$ such that
\smallskip

1. $E_T(\Lambda_1)E_T(\Lambda_2)=E_T(\Lambda_1\cap\Lambda_2)$,  
$E_T(\emptyset)=0$,  $E_T(\sigma(T))=I_{\cH}$.
\smallskip

2. $\|E_T(\Lambda)\|\leq C$ with $C$ independent of $\Lambda\in\cB$.
\smallskip

3. $E_T(\Lambda)\dom(T)\subseteq\dom(T)$,
$T E_T(\Lambda)f=E_T(\Lambda) Tf$, $f\in \dom(T)$, $\Lambda\in\cB$.
\smallskip

4. $E_T(\Lambda)\mathcal{H}\subseteq \dom(T)$ for $\Lambda$ bounded, 
$\sigma(T|_{E_T(\Lambda)\mathcal{H}\cap
\dom(T)})\subseteq\ol\Lambda$.
\smallskip

A spectral operator $T$ is a {\em spectral operator of scalar type} if
\begin{align*}
\dom(T)&=\bigg\{g\in\mathcal{H}\,\bigg|\,
\slimes_{n\uparrow\infty}\int_{\{\lambda\in\bbC | |\lambda|\leq n\}}
\lambda\,d(E_T(\lambda)g) \text{ exists in $\cH$}\bigg\} \\
Tf&=\slimes_{n\uparrow\infty}\int_{\{\lambda\in\bbC | |\lambda|\leq
n\}}
\lambda\,d(E_T(\lambda)f),
\quad f\in\dom(T).
\end{align*}

McGarvey was able to apply his results only to operators of the form
$-{d^2}/{dx^2} +p(x){d}/{dx} + q(x)$, $x\in\mathbb{R}$, 
with $\pi$-periodic functions $p$ and $q$ under the restriction
$\Im\bigg(\int_0^\pi dx \, p(x)\bigg) \neq0$. 
The spectra of such operators outside a sufficiently large disc
are composed of some separated ovals,
permitting McGarvey to prove that these operators are in some
sense  {\em asymptotically} spectral operators (\cite{Mc62}, Part II).
Such results ignore the existence of local spectral singularities and, at
any rate, are not applicable to Hill operators $H$.

Meiman \cite{Me77} noted that zeros of
$\Delta_+^\bullet(\lambda)$  are integrable singularities for the functions
\eqref{111} with $\Delta_+(z)=\cos(t)$, but this is generally incorrect
for $t\in\{0,\pi,2\pi\}$. Another claim of the same paper is  that ``the
crude asymptotic estimates for $\phi(\lambda,\pi)$,
$\Delta_+^\bullet(\lambda)$ and  the Floquet solutions of $H$ imply that
the Fourier integral theory remains essentially valid  also for expansions
in eigenfunctions of a non-self-adjoint Schr\"odinger (=Hill)  operator
with a complex-valued potential.'' However, such a vague statement,
without explanations as to what type of convergence and what type of
function space is meant, is unsatisfactory, especially, taking  into
account Theorems \ref{t2.1}--\ref{t2.3} stated below.

Finally, we mention that Veliev \cite{Ve83} erroneously concluded that the
spectra of Hill operators always have non-intersecting analytic arcs which
invalidates some his results concerning spectral  singularities and
spectral expansions.

Given this incomplete and, in part, quite confusing state of affairs on
the question of whether or not a Hill operator is a spectral operator
of scalar type, after more than 40 years since 
the problem first arose, we present its solution in terms of certain
functions related to  equation \eqref{20}.  For detailed proofs of Theorems
\ref{t2.1}--\ref{t2.3} we refer to \cite{GT06}.

\begin{theorem} \lb{t2.1}
A Hill operator $H$ is a spectral operator of scalar
type if  and only if the estimates
\begin{equation}\label{11}
\bigg|
\frac{\phi(\lambda,\pi)}{\Delta_+^{\bullet}(\lambda)}\bigg|\leq C,\quad
\left|
\frac{\theta'(\lambda,\pi)}{(|\lambda|+1)\Delta_+^{\bullet}(\lambda)}
\right|\leq C,\quad \bigg|
\frac{\Delta_-(\lambda)}{(\sqrt{|{\lambda}|}+1)
\Delta_+^{\bullet}(\lambda)}\bigg|\leq C
\end{equation}
hold for all  $\lambda\in\sigma(H)$, with $C$ a finite positive
constant independent of $\lambda\in\sigma(H)$.

If the conditions \eqref{11} are satisfied, then the functions \eqref{111}
are analytic for $z$ in an open neighborhood of $\sigma(H)$.
\end{theorem}

Some disadvantage of the above conditions lies in the fact that the
numerators of the fractions  in \eqref{11} are not independent functions.
To remedy this, we mention that Sansuc and Tkachenko \cite{ST96} gave
a parametrization of Hill operators using the functional parameters
$\{\phi(z,\pi),\Delta_+(z),\Delta_-(z)\}_{z\in\bbC}$. The following 
criterion is stated in terms of these parameters.

\begin{theorem} \lb{t2.2}
A Hill operator $H$ is a spectral operator of scalar type if and only if
the following conditions $(i)$ and $(ii)$ are satisfied: \\
$(i)$ The function
\begin{equation}\label{12}
\frac{\Delta_+(z)^2-1-\Delta_-(z)^2}
{\phi(z,\pi)\Delta_+^{\bullet}(z)}\quad
\end{equation}
is analytic for $z$ in an open neighborhood of $\sigma(H)$.\\
$(i)$ The inequalities
\begin{equation}\label{13}
\left|
\frac{\phi(\lambda,\pi)}{\Delta_+^{\bullet}(\lambda)}\right|\leq C,\quad
\bigg| \frac{\Delta_-(\lambda)}{(\sqrt{|{\lambda}|}+1)
\Delta_+^{\bullet}(\lambda)}\bigg|\leq C, \quad
\lambda\in\sigma(H),
\end{equation}
are satisfied with $C$ a finite positive constant independent of
$\lambda\in\sigma(H)$.
\end{theorem}

If both conditions \eqref{12} and \eqref{13} are satisfied, and a point
$\lambda_0\in\sigma(H)$ satisfies $\Delta_+^\bullet(\lambda_0)=0$, then 
$\Delta_+(\lambda_0)^2-1=\Delta_-(\lambda_0)=\Delta_+^\bullet(\lambda_0)
=0$, $\Delta_+^{\bullet\bullet}(\lambda_0)\neq 0$, 
implying that the spectrum of a Hill operator, which is a spectral operator
of scalar type, is formed by a system of countably many, 
non-intersecting, analytic arcs. The latter may degenerate into finitely
many simple analytic arcs and a simple analytic semi-infinite arc, all
of which are non-intersecting. Both conditions in \eqref{13} are
independent; neither one of them implies the other. Moreover, analyticity
of the function \eqref{12} does not follow from the estimates \eqref{13}.

To prove Theorems \ref{t2.1} and \ref{t2.2} we use, similar to the paper by
Tkachenko \cite{Tk64}, the method of direct integral decompositions
following Gel'fand \cite{Ge50}, connecting the Hill operator $H$ with
the family of densely defined, closed, linear operators $H(t)$, $t\in [0,
2\pi]$ in $L^2([0,\pi])$ defined by the differential expression
$-d^2/dx^2+q(x)$ restricted to $x\in [0,\pi]$ and the $t$-dependent
boundary conditions $y(\pi)=e^{it}y(0)$, $y'(\pi)=e^{it}y'(0)$. 
The spectrum of the operator $H(t)$, $t\in [0,2\pi]$, is given by
\begin{equation}
\sigma(H(t))=\{E_k(t)\}_{k\in\mathbb{N}_0}=\{z\in\mathbb{C}
\,|\, \Delta_+(z)=\cos (t)\},
\end{equation}
and the spectrum of $H$ is then given by
\begin{equation}
\sigma(H)=\bigcup\limits_{0\leq t\leq
\pi}\sigma(H(t)).
\end{equation}
For obvious reasons, the spectrum $\{E_k(0)\}_{k\in\bbN_0}$ of $H(0)$ and
$\{E_k(\pi)\}_{k\in\bbN_0}$ of $H(\pi)$ is called the periodic
and anti-periodic spectrum of $H$, respectively.

The following criterion involves the spectrum $\sigma(H)$, the periodic and
the anti-periodic spectra $\{\lambda_0^+, \lambda_k^\pm\}_{k\in\bbN}
=\{E_k(0), E_k(\pi)\}_{k\in\bbN_0}$ of $H$, the  Dirichlet spectrum
$\{\mu_k\}_{k\in\bbN}$, the set of critical points
$\{\delta_k\}_{k\in\bbN}$ of $\Delta_+$, and is connected with the
algebraic and geometric multiplicities of the eigenvalues in the sets
$\sigma(H(t))$, $t\in [0,\pi]$.

\begin{theorem} \lb{t2.3}
A Hill operator $H$ is a spectral operator of scalar type if and only if
the following conditions $(i)$--$(iii)$ are satisfied: \\
$(i)$  Every multiple point of either the periodic or anti-periodic
spectrum of $H$ is a point of its Dirichlet spectrum. \\
$(ii)$ For all $t\in[0,2\pi]$ and all $E_k(t)\in\sigma(H(t))$,
each root function (i.e., element of the algebraic eigenspace) of the
operator $H(t)$ associated with
$E_k(t)$ is an eigenfunction of $H(t)$. In particular, the geometric and
algebraic multiplicity of each eigenvalue $E_k(t)$ of $H(t)$ coincide. \\
$(iii)$ Let $
\mathcal{Q}=\{k\in\mathbb{N}\,|\, {\rm dist}(\delta_k,\sigma(H))\neq0\}$, 
then
\begin{equation} \label{15}
\displaystyle
\sup_{k\in\mathcal{Q}}\frac{|\lambda_k^+-\lambda_k^-|}
{{\rm dist}(\delta_k,\sigma(H))}<\infty,\quad
\sup_{k\in\mathcal{Q}}\frac{|\delta_k-\lambda_k^\pm|}
{{\rm dist}(\delta_k,\sigma(H))}<\infty.
\end{equation}
\end{theorem}

Our Theorems \ref{t2.1}--\ref{t2.3} give three equivalent criteria for the
eigenfunction system of the operators $H(t)$ to be a Riesz basis uniformly
with respect to $t\in[0,2\pi]$. The latter is equivalent for the operator
$H$ to be a spectral operator of scalar type.



\begin{thebibliography}{99}
%
\bi{DS88}  N.\ Dunford and J.\ T Schwartz, {\it Linear Operators, Part
III: Spectral Operators}, Wiley--Interscience, New York, 1988.
%
\bibitem {Ge50}   I.\ M.\ Gel'fand,
{\em Expansion in characteristic functions of an equation with periodic
coefficients}, Doklady Akad Nauk SSSR {\bf 73}, 1117-1120 (1950).
(Russian.)
%
\bi{GT06} F.\ Gesztesy and V.\ Tkachenko, {\it A criterion for Hill
operators to be spectral operators of scalar type}, in preparation.
%
\bi{GW95}  F.\ Gesztesy and R.\ Weikard, {\it Floquet theory revisited},
in {\it Differential Equations and Mathematical Physics}, I.\ Knowles
(ed.), International Press, Boston, 1995, pp.\ 67--84.
%
\bi{Mc62} D.\ McGarvey, {\it Operators commuting with translations
by one}. Part I, %
J. Math. Anal. Appl. {\bf 4},366-410 (1962),
%
Part II,
{\it ibid} {\bf 11}, 564-596 (1965); 
%
Part III, {\it ibid} 
{\bf 12}, 187--234 (1965).
%
\bi{Me77} N.\ N.\ Meiman, {\it The theory of one-dimensional Schr\"odinger
operators with a periodic potential}, J. Math. Phys. {\bf 18}, 834--848
(1977).
%
\bi{PT91a} L.\ A.\ Pastur and V.\ A.\ Tkachenko, {\it Geometry of the
spectrum of the one-dimensional Schr\"odinger equation with a periodic
complex-valued potential}, Math. Notes {\bf 50}, 1045--1050 (1991).
%
%
\bi{ST96} J.-J.\ Sansuc and V.\ Tkachenko, {\it Spectral
parametrization of non-selfadjoint Hill's operators}, J. Diff. Eq.
{\bf 125}, 366--384 (1996).
%
\bi{Se60} M.\ I.\ Serov, {\it Certain properties of the spectrum of a
non-selfadjoint differential operator of the second order}, Sov. Math.
Dokl. {\bf 1}, 190--192 (1960).
%
\bi{Ti50} E.\ C.\ Titchmarsh, {\it Eigenfunction problems with periodic
potentials}, Proc. Roy. Soc. London A {\bf 203}, 501--514 (1950).
%
\bi{Tk64} V.\ A.\ Tkachenko, {\it Spectral analysis of the one-dimensional
Schr\"odinger operator with periodic complex-valued potential}, Sov.
Math. Dokl. {\bf 5}, 413--415 (1964).
%
\bi{Tk96} V.\ A.\ Tkachenko, {\it Spectra of non-selfadjoint
Hill's operators and a class of Riemann surfaces}, Ann. Math. {\bf
143}, 181--231 (1996).
%
\bi{Ve83} O.\ A.\ Veliev, {\it Spectrum and spectral singularities of
differential operators with complex-valued periodic coefficients}, Diff.
Eqs. {\bf 19}, 983--989 (1983).
%
\end{thebibliography}
\end{document}